\documentclass[12pt]{article}
\usepackage{amsfonts}
\usepackage{amssymb}
\usepackage{amsmath}

\newcommand{\be}{\begin{equation}}
\newcommand{\ee}{\end{equation}}
\newcommand{\bea}{\begin{eqnarray*}}
\newcommand{\eea}{\end{eqnarray*}}
\newcommand{\ba}{\begin{array}}
\newcommand{\ea}{\end{array}}
\newcommand{\bi}{\begin{itemize}}
\newcommand{\ei}{\end{itemize}}
\newcommand{\bc}{\begin{center}}
\newcommand{\ec}{\end{center}}
\newcommand{\bfr}{\begin{flushright}}
\newcommand{\efr}{\end{flushright}}

\begin{document}

\title{Linear Continuous Functionals on {\bf FN}-Type Spaces}
\author{Sorin G. Gal \\
Department of Mathematics and Computer Sciences\\
University of Oradea, Romania\\
410087 Oradea, Romania\\
E-mail: galso@uoradea.ro}
\date{}
\maketitle

\begin{abstract}
By using the space of fuzzy numbers, in e.g. [5] have been
considered several complete metric spaces (called here {\bf
FN}-type spaces) endowed with addition and scalar multiplication,
such that the metrics have nice properties but the spaces are not
linear, i.e. are not groups with respect to addition and the
scalar multiplication is not, in general, distributive with
respect to usual scalar addition. This paper deals with the form
of linear continuous functionals defined on these spaces.
\end{abstract}


Subject Classification 2000 (AMS):46S40, 46A20, 46A45, 46E99.

\section{INTRODUCTION}

\quad\quad By using the space of fuzzy numbers, in [5] have been
considered several complete metric spaces endowed with addition,
and scalar multiplication, such that their metrics have nice
properties but the spaces are not linear, i.e. are not groups with
respect to addition and the scalar multiplication is not, in
general, distributive with respect to usual scalar addition. In
Section 2 we recall some properties of these spaces and introduce
the concept of abstract fuzzy-number-type ( shortly {\bf FN}-type)
space. Section 3 contains the main results of the paper and deals
with the form of linear and continuous functionals defined on the
{\bf FN}-type spaces in Section 2.

\section{PRELIMINARIES}

\quad\quad In this section we recall the main properties of the
space of fuzzy numbers and of some other spaces based on it, all
called as Fuzzy-Number-type (shortly {\bf FN}-type) spaces, which
have similar properties.

Given a set $X\neq \emptyset ,$ a fuzzy subset of $X$ is a mapping $%
u:X\rightarrow \left[ 0,1\right] $ and obviously any classical
subset $A$ of $X$ can be identified as a fuzzy subset of $X$
defined by $\chi _{A}:X\rightarrow \left[ 0,1\right] ,\chi
_{A}\left( x\right) =1$ if $x\in A,\chi _{A}\left( x\right) =0,$
if $x\in X\setminus A.$
If $u:X\rightarrow \left[ 0,1\right] $ is a fuzzy subset of $X,$ then for $%
x\in X,$ $u\left( x\right) $ is called the membership degree of
$x$ to $u$ (see e.g. $\left[ 14\right] $).

{\bf DEFINITION 2.1} (see e.g. $\left[ 4\right] $, $\left[
13\right] $)
The space of fuzzy numbers denoted by $\mathbf{%
R}_{\mathcal{F}}$ is defined as the class of fuzzy subsets of the real axis $%
\mathbf{R},$ i.e. of $u:\mathbf{R}\rightarrow \left[ 0,1\right] ,$ having
the following four properties:

$\left( i\right) $ $\forall u\in \mathbf{R}_{\mathcal{F}},$ $u$ is
normal, i.e. $\exists x_{u}\in \mathbf{R}$ with $u\left(
x_{u}\right) =1;$

$\left( ii\right) $ $\forall u\in \mathbf{R}_{\mathcal{F}},$ $u$ is a convex
fuzzy set, i.e.
\[
u\left( tx+\left( 1-t\right) y\right) \geq \min \left\{ u\left(
x\right)
,u\left( y\right) \right\} ,\forall t\in \left[ 0,1\right] ,x,y\in \mathbf{R}%
;
\]

$\left( iii\right) $ $\forall u\in \mathbf{R}_{\mathcal{F}},$ $u$
is upper-semi-continuous on $\mathbf{R};$

$\left( iv\right) $ $\overline{\left\{ x\in
\mathbf{R}\text{;}u\left( x\right) >0\right\} }$ is compact, where
$\overline{M}$ denotes the closure of $M.$

{\bf REMARKS.} 1) Obviously, we can consider that
$\mathbf{R}\subset \mathbf{R}_{\mathcal{F}},$ because any real
number $x_{0}\in \mathbf{R}$ can be identified with $\chi
_{\left\{ x_{0}\right\} },$ which satisfies the properties $\left(
i\right) -\left( iv\right) $ in Definition 2.1.

 2) For $0<r\leq 1$ and $u\in \mathbf{R}_{\mathcal{F}},$ let us denote by $%
\left[ u\right] ^{r}=\left\{ x\in \mathbf{R};u\left( x\right) \geq
r\right\} $ and $\left[ u\right] ^{0}=\overline{\left\{ x\in
\mathbf{R};u\left( x\right) >0\right\} },$ the so-called level
sets of $u.$ Then it is an immediate consequence of $\left(
i\right) -\left( iv\right) $ that $\left[
u\right] ^{r}$ represents a bounded closed (i.e. compact) subinterval of $%
\mathbf{R},$ denoted by $\left[ u\right] ^{r}=\left[ u_{-}\left(
r\right) ,u_{+}\left( r\right) \right] ,$ where $u_{-}\left(
r\right) \leq u_{+}\left( r\right) $ for all $r\in \left[
0,1\right] .$ Also, by e.g. [10], [13], $u_{-}\left( r\right) $ is
bounded nondecreasing on $\left[ 0,1\right] ,$ $u_{+}\left(
r\right) $ is bounded non-increasing on $\left[ 0,1\right] ,$ both
are left continuous on $\left( 0,1\right] $ and right continuous
at $r=0$, (from monotonicity both have right limit at each point
in $[0,1]$ ), $ u_{+}\left( 0\right) -u_{+}\left( r\right) \geq
0,u_{-}\left( 1\right) -u_{-}\left( r\right) \geq 0$, $u_{-}\left(
r\right)\le u_{+}\left( r\right), \forall r\in [0,1]$ and
$\mathbf{R}_{\mathcal{F}}$ can be embeded into the Banach space
$B={\overline C}[0,1]\times {\overline C}[0,1]$, by the mapping
$j(u)=(u_{-},u_{+}), \forall u\in \mathbf{R}_{\mathcal{F}}$, where
${\overline C}[0,1]$ denotes the Banach space of all real-valued
bounded functions $f:[0,1]\to \mathbb{R}$, which are left
continuous at each point in $(0,1]$, have right limit at each
point in $[0,1]$, $f$ is right continuous at $0$, endowed with the
uniform norm $||f||=\sup \{|f(x)| ; x\in [0,1]\}$ and the product
space $B$ is considered to be endowed with the norm
$||(f,g)||=\max \{||f||, ||g||\}$.

Also, it is important the following "converse" result.

{\bf THEOREM 2.2} (see e.g. [9] or [13, Lemma 1.1] ){\it If
$\{M_{r} ; r\in [0,1]\}$ is a family of closed subintervals of
real axis with the properties :

(i) $M_{r}\subset M_{s}, \forall r,s \in [0,1], s\le r$,

(ii) for any sequence $(r_{n})_{n\in \mathbb{N}}$ converging
increasingly to $r\in (0,1]$, we have $\bigcap
_{n=1}^{\infty}M_{r_{n}}=M_{r}$,

then there exists a unique $u\in \mathbf{R}_{\mathcal{F}}$, such
that $M_{r}=[u_{-}(r),u_{+}(r)], \forall r\in (0,1]$ and
$[u]^{0}\subset M_{0}$.}

{\bf DEFINITION 2.3} (see e.g. $\left[ 4\right] $, $\left[
13\right] $)
The addition and the product with real scalars in $%
\mathbf{R}_{\mathcal{F}}$ are defined by $\oplus :\mathbf{R}_{\mathcal{F}%
}\times \mathbf{R}_{\mathcal{F}}\rightarrow \mathbf{R}_{\mathcal{F}},$%
\[
\left( u\oplus v\right) \left( x\right) =\sup_{y+z=x}\min \left\{
u\left( y\right) ,v\left( z\right) \right\}
\]
and by $\odot :\mathbf{R}\times \mathbf{R}_{\mathcal{F}}\rightarrow \mathbf{R%
}_{\mathcal{F}},$%
\[
\left( \lambda \odot v\right) \left( x\right) =\left\{
\begin{array}{l}
u\left( \frac{x}{\lambda }\right)  \\
\widetilde{0}
\end{array}
\begin{array}{l}
\text{if }\lambda \neq 0 \\
\text{if }\lambda =0
\end{array}
\right. ,
\]
where $\widetilde{0}:\mathbf{R}\rightarrow \left[ 0,1\right] $ is $%
\widetilde{0}=\chi _{\left\{ 0\right\} }.$

Also, we can write $\left[ u\oplus v\right] ^{r}=\left[ u\right]
^{r}+\left[ v\right] ^{r},\left[ \lambda \odot v\right]
^{r}=\lambda \left[ v\right] ^{r},$ for all $r\in \left[
0,1\right] ,$ where $\left[ u\right] ^{r}+\left[ v\right] ^{r}$
means the usual sum of two intervals (as subsets of $\mathbf{R}$)
and $\lambda \left[ v\right] ^{r}$ means the usual product between
a real scalar and a subset of $\mathbf{R}.$

If we define $D:\mathbf{R}_{\mathcal{F}}\times \mathbf{R}_{\mathcal{F}%
}\rightarrow \mathbf{R}_{+}\cup \left\{ 0\right\} $ by
\[
D\left( u,v\right) =\sup_{r\in \left[ 0,1\right] }\max \left\{ \left|
u_{-}\left( r\right) -v_{-}\left( r\right) \right| ,\left| u_{+}\left(
r\right) -v_{+}\left( r\right) \right| \right\} ,
\]
where $\left[ u\right] ^{r}=\left[ u_{-}\left( r\right) ,u_{+}\left(
r\right) \right] ,$ $\left[ v\right] ^{r}=\left[ v_{-}\left( r\right)
,v_{+}\left( r\right) \right] ,$ then we have the following:

{\bf THEOREM 2.4} ( see e.g. $\left[ 10\right] $, $\left[
13\right] $). {\it $\left( \mathbf{R}_{\mathcal{F}},D\right) $ is
a complete metric space and in addition, $D$ has the following
three properties:

$\left( i\right) $ $D\left( u\oplus w, v\oplus w\right) =D\left(
u,v\right) ,$ for all $u,v,w\in \mathbf{R}_{\mathcal{F}};$

$\left( ii\right) $ $D\left( k\odot u,k\odot v\right) =\left| k\right|
D\left( u,v\right) ,$ for all $u,v\in \mathbf{R}_{\mathcal{F}},k\in \mathbf{R%
};$

$\left( iii\right) $ $D\left( u\oplus v,w\oplus e\right) \leq D\left(
u,w\right) +D\left( v,e\right) ,$ for all $u,v,w,e\in \mathbf{R}_{\mathcal{F}%
}.$}

Also, the following result is known:

{\bf THEOREM 2.5} (see e.g. $\left[ 1\right] ,\left[ 4\right] $).
{\it $\left( i\right) $ $u\oplus v=v\oplus u,u\oplus \left(
v\oplus w\right) =\left( u\oplus v\right) \oplus w;$

$\left( ii\right) $ If we denote $\widetilde{0}=\chi _{\left\{
0\right\} },$ then $u\oplus \widetilde{0}=\widetilde{0}\oplus
u=u,$ for any $u\in \mathbf{R}_{\mathcal{F}};$

$\left( iii\right) $ With respect to $\widetilde{0},$ none of $u\in \mathbf{R%
}_{\mathcal{F}}\setminus \mathbf{R}$ has an opposite element (regarding $%
\oplus $) in $\mathbf{R}_{\mathcal{F}};$

$\left( iv\right) $ For any $a,b\in \mathbf{R}$ with $a,b\geq 0$
or $a,b\leq 0$ and any $u\in \mathbf{R}_{\mathcal{F}},$ we have
\[
\left( a+b\right) \odot u=a\odot u\oplus b\odot u.
\]
For general $a,b\in \mathbf{R},$ the above property does not hold.

$\left( v\right) $ $\lambda \odot \left( u\oplus v\right) =\lambda \odot
u\oplus \lambda \odot v,$ for all $\lambda \in \mathbf{R},u,v\in \mathbf{R}_{%
\mathcal{F}};$

$\left( vi\right) $ $\lambda \odot \left( \mu \odot u\right) =\left( \lambda
\mu \right) \odot u,$ for all $\lambda ,\mu \in \mathbf{R},u\in \mathbf{R}_{%
\mathcal{F}};$

$\left( vii\right) $ If we denote $\left\| u\right\| _{\mathcal{F}}=D\left(
u,\widetilde{0}\right) ,u\in \mathbf{R}_{\mathcal{F}},$ then $\left\|
u\right\| _{\mathcal{F}}$ has the properties of an usual norm on $\mathbf{R}_{%
\mathcal{F}},$ i.e. $\left\| u\right\| _{\mathcal{F}}=0$ iff $u=\widetilde{0}%
,\left\| \lambda \odot u\right\| _{\mathcal{F}}=\left| \lambda \right|
\left\| u\right\| _{\mathcal{F}},\left\| u+v\right\| _{\mathcal{F}}\leq
\left\| u\right\| _{\mathcal{F}}+\left\| v\right\| _{\mathcal{F}},\left|
\left\| u\right\| _{\mathcal{F}}-\left\| v\right\| _{\mathcal{F}}\right|
\leq D\left( u,v\right) ;$

$\left( viii\right) D\left( \alpha \odot u,\beta \odot u\right)
=\left| \alpha -\beta \right| D\left( \widetilde{0},u\right) ,$
for all $\alpha ,\beta \geq 0,u\in \mathbf{R}_{\mathcal{F}}.$ If
$\alpha ,\beta \leq 0$ then the equality is also valid. If $\alpha
$ and $\beta $ are of opposite signs, then the equality is not
valid.}

{\bf REMARKS.} 1) Theorem 2.5 shows that $\left(
\mathbf{R}_{\mathcal{F}},\oplus ,\odot \right) $ is not a linear
space over $\mathbf{R}$ and consequently $\left(
\mathbf{R}_{\mathcal{F}},\left\| u\right\| _{\mathcal{F}}\right) $
cannot be a normed space.

2) On $\mathbf{R}_{\mathcal{F}}$ we also can define a substraction
$\ominus ,$ called $H-$ difference (see e.g. $\left[ 3\right] $)
as follows: $u\ominus v$
has sense if there exists $w\in \mathbf{R}_{\mathcal{F}}$ such that $%
u=v\oplus w.$ Obviously, $u\ominus v$ does not exist for all
$u,v\in \mathbf{R}_{\mathcal{F}}$ (for example,
$\widetilde{0}\ominus v$ does not exists if $v\neq
\widetilde{0}$).

In what follows, we define some usual spaces of
fuzzy-number-valued
functions, which have similar properties to $\left( \mathbf{R}_{\mathcal{F}%
},D\right) .$

Denote $C\left( \left[ a,b\right] ;\mathbf{R}_{\mathcal{F}}\right) =\left\{
f:\left[ a,b\right] \rightarrow \mathbf{R}_{\mathcal{F}};\text{ f is
continuous on }\left[ a,b\right] \right\} ,$ endowed with the metric $%
D^{*}\left( f,g\right) =\sup \left\{ D\left( f\left( x\right) ,g\left(
x\right) \right) ;x\in \left[ a,b\right] \right\} .$ Because $\left( \mathbf{%
R}_{\mathcal{F}},D\right) $ is a complete metric space, by
standard technique we obtain that $\left( C\left( \left[
a,b\right] ;\mathbf{R}_{\mathcal{F}}\right) ,D^{*}\right) $ is a
complete metric space. Also, if we define $\left( f\oplus g\right)
\left( x\right) =f\left( x\right) \oplus g\left( x\right) ,\left(
\lambda \odot f\right) \left( x\right) =\lambda \odot f\left(
x\right) $ (for simplicity,
the addition and scalar multiplication in $C\left( \left[ a,b\right] ;%
\mathbf{R}_{\mathcal{F}}\right) $ are denoted as in $\mathbf{R}_{\mathcal{F}%
} $), also $\widetilde{0}:\left[ a,b\right] \rightarrow
\mathbf{R}_{\mathcal{F}},\widetilde{0}\left(t\right)=\widetilde{0}_{\mathbf{R}_{\mathcal{F}}},$
for all $t\in \left[ a,b\right] ,$%
\[
\left\| f\right\| _{\mathcal{F}}=\sup \left\{ D\left(
\widetilde{0},f\left( x\right) \right) ;x\in \left[ a,b\right]
\right\} ,
\]
then we easily obtain the following properties.

{\bf THEOREM 2.6} (see $\left[ 5\right] $) {\it $\left( i\right) $
$f\oplus g=g\oplus f,\left( f\oplus g\right) \oplus h=f\oplus
\left( g\oplus h\right) ;$

$\left( ii\right) $ $f\oplus \widetilde{0}=\widetilde{0}\oplus f,$ for any $%
f\in C\left( \left[ a,b\right] ;\mathbf{R}_{\mathcal{F}}\right) ;$

$\left( iii\right) $ With respect to $\widetilde{0}$ in $C\left(
\left[ a,b\right] ;\mathbf{R}_{\mathcal{F}}\right) ,$ any $f\in
C\left( \left[ a,b\right] ;\mathbf{R}_{\mathcal{F}}\right) $ with
$f\left( \left[ a,b\right] \right) \cap
\mathbf{R}_{\mathcal{F}}\neq \emptyset $ has no an
opposite member (regarding $\oplus $) in $C\left( \left[ a,b\right] ;\mathbf{%
R}_{\mathcal{F}}\right) ;$

$\left( iv\right) $ for all $\lambda ,\mu \in \mathbf{R}$ with $\lambda ,\mu
\geq 0$ or $\lambda ,\mu \leq 0$ and for any $f\in C\left( \left[ a,b\right]
;\mathbf{R}_{\mathcal{F}}\right) ,$%
\[
\left( \lambda +\mu \right) \odot f=\left( \lambda \odot f\right) \oplus
\left( \mu \odot f\right) ;
\]
For general $\lambda ,\mu \in \mathbf{R},$ this property does not hold.

$\left( v\right) $ $\lambda \odot \left( f\oplus g\right) =\lambda
\odot f\oplus \lambda \odot g,\lambda \odot \left( \mu \odot
f\right) =\left(
\lambda \mu \right) \odot f,$ for any $f,g\in C\left( \left[ a,b\right] ;%
\mathbf{R}_{\mathcal{F}}\right) ,\lambda ,\mu \in \mathbf{R};$

$\left( vi\right) $ $\left\| f\right\| _{\mathcal{F}}=0$ iff $f=\widetilde{0}%
,$ $\left\| \lambda \odot f\right\| _{\mathcal{F}}=\left| \lambda \right|
\left\| f\right\| _{\mathcal{F}},\left\| f\oplus g\right\| _{\mathcal{F}%
}\leq \left\| f\right\| _{\mathcal{F}}+\left\| g\right\| _{\mathcal{F}%
},\left|  \left\| f\right\| _{\mathcal{F}}-\left\| g\right\| _{\mathcal{F}%
}\right| \leq D^{*}\left( f,g\right) ,$ for any $f,g\in C\left( \left[
a,b\right] ;\mathbf{R}_{\mathcal{F}}\right) ,\lambda \in \mathbf{R};$

$\left( vii\right) $ $D^{*}\left( \lambda \odot f,\mu \odot
f\right) =\left| \lambda -\mu \right| D^{*}\left(
\widetilde{0},f\right) ,$for any $f\in C\left( \left[ a,b\right]
;\mathbf{R}_{\mathcal{F}}\right) ,\lambda \mu \geq 0;$

$\left( viii\right) $%
\[
D^{*}\left( f\oplus h,g\oplus h\right) =D^{*}\left( f,g\right) ,
\]
\[
D^{*}\left( \lambda \odot f,\lambda \odot g\right) =\left| \lambda \right|
D^{*}\left( f,g\right) ,
\]
\[
D^{*}\left( f\oplus g,h\oplus e\right) \leq D^{*}\left( f,h\right)
+D^{*}\left( g,e\right) ,
\]
for any $f,g,h,e\in C\left( \left[ a,b\right] ;\mathbf{R}_{\mathcal{F}%
}\right) ,\lambda \in \mathbf{R}.$}

{\bf REMARK.}
It is easy to show that if $f,g\in C\left( \left[ a,b\right] ;\mathbf{R}_{%
\mathcal{F}}\right) ,$ then $F:\left[ a,b\right] \rightarrow
\mathbf{R},$ defined by $F\left( x\right) =D\left( f\left(
x\right) ,g\left( x\right) \right) $ is continuous on $\left[
a,b\right] .$

Now, for $1\leq p<\infty ,$ let us define
\[
L^{p}\left( \left[ a,b\right] ;\mathbf{R}_{\mathcal{F}}\right) =\left\{
\begin{array}{c}
f\text{ is strongly measurable on }\left[ a,b\right] \text{ and } \\
\left( L\right) \int\limits_{a}^{b}\left( D\left( \widetilde{0},f\left(
x\right) \right) \right) ^{p}dx<+\infty
\end{array}
\right\} ,
\]
where according to e.g. $\left[ 8\right] $, $f$ is called strongly
measurable if, for each $x\in \left[ a,b\right] ,$ $f_{-}\left(
x\right) \left( r\right) $ and $f_{+}\left( x\right) \left(
r\right) $ are Lebesgue measurable as functions of $r\in \left[
0,1\right] $ (here $\left[ f\left( x\right) \right] ^{r}=\left[
f_{-}\left( x\right) \left( r\right) ,f_{+}\left( x\right) \left(
r\right) \right] $ denotes the $r-$level set of $f\left( x\right)
\in \mathbf{R}_{\mathcal{F}}$). The following result shows that
$L^{p}\left( \left[ a,b\right] ;\mathbf{R}_{\mathcal{F}}\right) $
is well defined.

{\bf THEOREM 2.7} (see $\left[ 5\right] $)
{\it $\left( i\right) $ If $f:\left[ a,b\right] \rightarrow \mathbf{R}_{\mathcal{F%
}}$ is strongly measurable then $F:\left[ a,b\right] \rightarrow \mathbf{R}%
_{+}$ defined by $F\left( x\right) =D\left( \widetilde{0},f\left(
x\right) \right) $ is Lebesgue measurable on $\left[ a,b\right] ;$

$\left( ii\right) $ For any $f,g\in L^{p}\left( \left[ a,b\right] ;\mathbf{R}%
_{\mathcal{F}}\right) ,F\left( x\right) =D\left( f\left( x\right)
,g\left( x\right) \right) $ is Lebesgue measurable and
$L^{p}$-integrable on $\left[ a,b\right] .$ Moreover, if we define
\[
D_{p}\left( f,g\right) =\left\{ \left( L\right) \int\limits_{a}^{b}\left[
D\left( f\left( x\right) ,g\left( x\right) \right) \right] ^{p}dx\right\} ^{%
\frac{1}{p}},
\]
then $\left( L^{p}\left( \left[ a,b\right]
;\mathbf{R}_{\mathcal{F}}\right) ,D_{p}\right) $ is a complete
metric space (where $f=g$ means $f\left( x\right) =g\left(
x\right) ,$ a.e. $x\in \left[ a,b\right] $) and, in addition,
$D_{p} $ satisfies the following properties:
\[
D_{p}\left( f\oplus h,g\oplus h\right) =D_{p}\left( f,g\right) ,
\]
\[
D_{p}\left( \lambda \odot f,\lambda \odot g\right) =\left| \lambda \right|
D_{p}\left( f,g\right) ,
\]
\[
D_{p}\left( f\oplus g,h\oplus e\right) \leq D_{p}\left( f,h\right)
+D_{p}\left( g,e\right) ,
\]
for any $f,g,h,e\in L^{p}\left( \left[ a,b\right] ;\mathbf{R}_{\mathcal{F}%
}\right) ,\lambda \in \mathbf{R}.$}

Other spaces with properties similar to those of $\left(
\mathbf{R}_{\mathcal{F}},D\right) $ can be constructed as follows
(see $\left[ 5\right] $).

For $p\geq 1,$ let us define
\[
l_{\mathbf{R}_{\mathcal{F}}}^{p}=\left\{ x=\left( x_{n}\right) _{n};x_{n}\in
\mathbf{R}_{\mathcal{F}},\forall n\in \Bbb{N}\text{ and }\sum\limits_{n=1}^{%
\infty }\left\| x_{n}\right\| _{\mathbf{R}_{\mathcal{F}}}^{p}<+\infty
\right\} ,
\]
endowed with the metric
\[
\rho _{p}\left( x,y\right) =\left\{ \sum\limits_{n=1}^{\infty
}[D\left( x_{n},y_{n}\right)]^{p} \right\} ^{\frac{1}{p}},\forall
x=\left( x_{n}\right) _{n},y=\left( y_{n}\right) _{n}\in
l_{\mathbf{R}_{\mathcal{F}}}^{p}.
\]
By $D\left( x_{n},y_{n}\right) =D\left( x_{n}\oplus \widetilde{0},\widetilde{%
0}\oplus y_{n}\right) \leq D\left( x_{n},\widetilde{0}\right) +D\left(
\widetilde{0},y_{n}\right) =\left\| x_{n}\right\| _{\mathcal{F}}+\left\|
y_{n}\right\| _{\mathcal{F}},$ we easily get (by Minkowski's inequality if $%
p>1$) $\rho _{p}\left( x,y\right) <+\infty .$ Also, it easily follows that $%
\rho _{p}\left( x,y\right) $ is a metric with similar properties
to $D$ (see Theorem 2.4 and Theorem 2.5,(vii),(viii)).

Because $\left( \mathbf{R}_{\mathcal{F}},D\right) $ is a complete metric
space, by the standard technique, we easily get that $\left( l_{\mathbf{R}_{%
\mathcal{F}}}^{p},\rho _{p}\right) $ is also a complete metric space.

Let us denote by
\[
m_{\mathbf{R}_{\mathcal{F}}}=\left\{ x=\left( x_{n}\right) _{n};x_{n}\in
\mathbf{R}_{\mathcal{F}},\forall n\in \Bbb{N}\text{ and }\exists M>0\text{
such that }\left\| x_{n}\right\| _{\mathcal{F}}\leq M,\forall n\in \Bbb{N}%
\right\} ,
\]
endowed with the metric
\[
\mu \left( x,y\right) =\sup \left\{ D\left( x_{n},y_{n}\right) ,\forall n\in
\Bbb{N}\right\} .
\]
We easily get that $\left( m_{\mathbf{R}_{\mathcal{F}}},\mu
\right) $ is a complete metric space and, in addition, $\mu $ has
similar properties to $D$ (see Theorem 2.4 and Theorem
2.5,(vii),(viii)). Similarly, if we denote
\[
c_{\mathbf{R}_{\mathcal{F}}}=\left\{ x=\left( x_{n}\right) _{n};x_{n}\in
\mathbf{R}_{\mathcal{F}},\forall n\in \Bbb{N}\text{ and }\exists a\in
\mathbf{R}_{\mathcal{F}}\text{ such that }D\left( x_{n},a\right) \stackrel{%
n\rightarrow \infty }{\longrightarrow }0\right\}
\]
and
\[
c_{\mathbf{R}_{\mathcal{F}}}^{\widetilde{0}}=\left\{ x=\left( x_{n}\right)
_{n};x_{n}\in \mathbf{R}_{\mathcal{F}},\forall n\in \Bbb{N}\text{, such that
}D\left( x_{n},\widetilde{0}\right) \stackrel{n\rightarrow \infty }{%
\longrightarrow }0\right\} ,
\]
since $\left( \mathbf{R}_{\mathcal{F}},D\right) $ is complete, by
standard technique, it follows that $\left(
c_{\mathbf{R}_{\mathcal{F}}},\mu \right) $ and $\left(
c_{\mathbf{R}_{\mathcal{F}}}^{\widetilde{0}},\mu \right) $ are
complete metric spaces.

{\bf REMARK}. Let $\left( X,\oplus ,\odot ,d\right) $ be represent
any space from $\left( \mathbf{R}_{\mathcal{F}},D\right) ,\left(
l_{\mathbf{R}_{ \mathcal{F}}}^{p},\rho _{p}\right) ,\left(
m_{\mathbf{R}_{\mathcal{F}}},\mu
\right) ,\left( c_{\mathbf{R}_{\mathcal{F}}},\mu \right) ,\left( c_{\mathbf{R%
}_{\mathcal{F}}}^{\widetilde{0}},\mu \right) ,\left( L^{p}\left(
\left[ a,b\right] ;\mathbf{R}_{\mathcal{F}}\right) ,D_{p}\right)
,1\leq p<\infty $, $\ \left( C\left( \left[ a,b\right]
;\mathbf{R}_{\mathcal{F}}\right) ,D^{*}\right) ,$ or any finite
cartesian product of them. The properties in Theorems 2.4, 2.5,
2.6, 2.7, suggest us in a natural way the following concept of
abstract space.

{\bf DEFINITION 2.8} We say that $\left( X,\oplus, \odot ,d\right)
$ is a fuzzy-number type space (shortly {\bf FN}-type space), if
the following properties are satisfied :

(i) $(X,d)$ is a metric space (complete or not) and $d$ has the
properties in Theorem 2.4, (i)-(iii) (where
$\mathbf{R}_{\mathcal{F}}$  is replaced by $X$ and $D$ by $d$ );

(ii) The operations $\oplus, \odot$ on $X$ have the properties in
Theorem 2.5, (i),(iv),(v),(vi) (where $\mathbf{R}_{\mathcal{F}}$
is replaced by $X$ ) ;

(iii) There exists a neutral element $\widetilde{0}\in X$, i.e.
$u\oplus \widetilde{0}=\widetilde{0}\oplus u=u,$ for any $u\in X$
and a linear subspace $Y\subset X$ (with respect to $\oplus$ and
$\odot$), non-dense in $X$, such that with respect to
$\widetilde{0},$ none of
$u\in X\setminus Y$ has an opposite member (regarding $%
\oplus $) in $X$.

{\bf REMARK.} A {\bf FN}-type space obviously is a more general
structure that that of Banach or Fr\'echet space, because it is
not a linear space. However, due to the nice properties of the
metric $d$, very many results (especially those of quantitative
kind) valid for Banach (or Fr\'echet) spaces, can be extended to
this case too. For example, the theories of almost periodic and
almost automorphic functions with values ${\bf FN}$-type spaces
were developed in [2] and [7], respectively. Also, the ${\bf
FN}$-type spaces have recent applications to the study of fuzzy
differential equations, which model the real world's problems
governed by imprecision due to uncertainty or vagueness rather
than randomness. In this sense, we mention for example [5] and
[6], where basic elements of the theory of semigroups of operators
on ${\bf FN}$-type spaces with applications in solving fuzzy
partial differential equations are considered. Of course that the
theory of semigroups of operators requires basic elements of
operator theory on these spaces, as for example, the followings.

{\bf DEFINITION 2.9} ($\left[ 5\right] $) $A:X\Bbb{\rightarrow
}\mathbf{R}$ is a linear functional if
\[
\left\{
\begin{array}{l}
A\left( x\oplus y\right) =A\left( x\right) +A\left( y\right) , \\
A\left( \lambda \odot x\right) =\lambda A\left( x\right) ,
\end{array}
\right.
\]
for all $x,y\in X$,$\lambda \in \mathbf{R}$.

{\bf REMARK.} If $A:X\Bbb{\rightarrow }\mathbf{R}$ is linear and
continuous at $\widetilde{0}\in X$, then this does not imply the
continuity of $A$ at each $x\in X$, because we cannot write
$x_{0}=(x_{0}\ominus x)\oplus x$ ,in general,$ $(the difference
$x_{0}\ominus x$ does not always exist).

However, we can prove the following theorem.

{\bf THEOREM 2.10} ($\left[ 5\right] $) {\it If
$A:X\Bbb{\rightarrow }\mathbf{R}$ is linear, then it is continuous
at $\widetilde{0}\in X$, if and only if there exists $M>0$ such
that
\[
\left| A\left( x\right) \right| \leq M\left\| x\right\| _{\mathcal{F}%
},\forall x\in X\text{,}
\]
where $\left\| x\right\| _{\mathcal{F}}=d\left(
\widetilde{0},x\right) .$}

Now, for $A:X\Bbb{\rightarrow }\mathbf{R}$ linear and continuous at $%
\widetilde{0},$ let us denote by
\[
\mathcal{M}_{A}=\left\{ M>0;\left| A\left( x\right) \right| \leq
M\left\| x\right\| _{\mathcal{F}},\forall x\in X\right\} ,
\]

Also, denote $\left\| \left| A\right| \right\| _{\mathcal{F}}=%
\inf \mathcal{M}_{A}.$

We have

{\bf THEOREM 2.11} (see $\left[ 5\right] $) {\it If
$A:X\Bbb{\rightarrow }\mathbf{R}$ is linear and continuous at
$\widetilde{0},$ then
\[
\left| A\left( x\right) \right| \leq \left| \left\| A\right\| \right| _{%
\mathcal{F}}\left\| x\right\| _{\mathcal{F}}
\]
for all $x\in X$ and
\[
\left| \left\| A\right\| \right| _{\mathcal{F}}=\sup \left\{
\left| A\left( x\right) \right| ;x\in X\text{,}\left\| x\right\|
_{\mathcal{F}}\leq 1\right\} .
\]}

{\bf COROLLARY 2.12} (see $\left[ 5\right] $) {\it If $A:$
$X\Bbb{\rightarrow }$ $\mathbf{R}$ is additive (i.e. $A\left(
x\oplus y\right) =A\left( x\right) +A\left( y\right) $), positive
homogeneous (i.e. $A\left( \lambda \odot y\right) =\lambda A\left(
x\right) ,\forall \lambda \geq 0$) and continuous at
$\widetilde{0},$ then
\[
\left| A\left( x\right) \right| \leq \left| \left\| A\right\| \right| _{%
\mathcal{F}}\left\| x\right\| _{\mathcal{F}},\forall x\in
X\text{.}
\]}

Also, the following uniform boundedness principle holds.

{\bf THEOREM 2.13} (see [5]) \textit{Let $\left( X,\oplus, \cdot,
d\right) $ be a ${\bf FN}$-type space and $\mathbf{L}\left(
X\right) $ be any from the spaces
\[
\mathcal{L}^{+}\left( X\right) =\left\{ A\in
\mathcal{L}_{0}^{+}\left( X\right) ;A\text{ is continuous at each
}x\in X\right\} ,
\]
\[
\mathcal{L}\left( X\right) =\left\{ A\in \mathcal{L}_{0}\left( X\right) ;A%
\text{ is continuous at each }x\in X\right\} ,
\]
where
\[
\mathcal{L}_{0}^{+}\left( X\right) =\left\{ A:X\Bbb{\rightarrow
}X;A\text{ is additive, positive homogeneous and continuous at
}\widetilde{0}\right\} ,
\]
$\mathcal{L}_{0}\left( X\right) =\left\{ A:X\Bbb{\rightarrow
}X;A\text{ is linear and continuous at }\widetilde{0}\right\} .$}

\textit{If $A_{j}\in \mathbf{L}\left( X\right) ,j\in J,$ is
pointwise bounded, i.e. for any $x\in X$, $\left\| A_{j}\left( x\right) \right\| _{%
\mathcal{F}}=d\left( A_{j}\left( x\right) ,\widetilde{0}\right)
\leq M_{x},$ for all $j\in J,$ then there exists $M>0$ such that
\[
\left| \left\| A_{j}\right\| \right| _{\mathcal{F}}\leq M,\forall j\in J%
\text{,}
\]
(i.e. $\left( A_{j}\right) _{j}$ is uniformly bounded).}

{\bf REMARK.} It is worth to note that not all the results in
operator theory on Banach spaces can be extended to ${\bf
FN}$-type spaces (see [5]).

\section{FORMS OF THE LINEAR CONTINUOUS FUNCTIONALS}
\quad \quad
In this section we deal with the form of linear and
continuous functionals defined on the ${\bf FN}$-type spaces
introduced by Section 2. The first main result is the following.

{\bf THEOREM 3.1 } {\it $ x^{*}: \mathbf{R}_{\mathcal{F}}\to
\mathbf{R}$ is a linear continuous functional on
$\mathbf{R}_{\mathcal{F}}$, if and only if there exists a linear
functional $L : \overline{C}[0,1]\to \mathbf{R}$, such that
$$x^{*}(x)=L(x_{-}+x_{+}), \forall x\in
\mathbf{R}_{\mathcal{F}},$$ where $x_{-}$ and $x_{+}$ are the
functions given by the formula $\left[ x\right] ^{r}=\left[
x_{-}\left( r\right) ,x_{+}\left( r\right) \right] ,$ (see Remark
2 after Definition 2.1) and $L_{|IC[0,1]}$ is continuous with
respect to the uniform convergence on $IC[0,1]=\{f\in
\overline{C}[0,1] : f\text {is increasing on }  [0,1]\}$. (Here
$L_{|A}$ denotes the restriction of $L$ to $A$ ).}

{\it Proof.} Let $L$ be as in the statement and define
$x^{*}(x)=L(x_{-}+x_{+}), \forall x\in \mathbf{R}_{\mathcal{F}}.$
First we prove that $x^{*}$ is linear.

Let $x, y\in \mathbf{R}_{\mathcal{F}}$. From the obvious relations
$(x\oplus y)_{-}=x_{-}+y_{-}, (x\oplus y)_{+}=x_{+}+y_{+}$ and the
linearity of $L$, we immediately get $x^{*}(x\oplus
y)=x^{*}(x)+x^{*}(y)$.

Let $\alpha\in \mathbf{R}$ and $x\in \mathbf{R}_{\mathcal{F}}$. If
$\alpha\ge 0$ then by the obvious relations $(\alpha\odot
x)_{-}=\alpha (x)_{-}, (\alpha\odot x)_{+}=\alpha (x)_{+}$ and the
linearity of $L$, we easily obtain $x^{*}(\alpha\odot x)=\alpha
x^{*}(x)$. If $\alpha < 0$ then by the relations $(\alpha\odot
x)_{-}=\alpha (x)_{+}, (\alpha\odot x)_{+}=\alpha (x)_{-}$ and the
linearity of $L$, we again arrive at the same conclusion
$x^{*}(\alpha\odot x)=\alpha x^{*}(x)$.

Next we prove that $x^{*}$ is continuous. For that let $x_{n},
x\in \mathbf{R}_{\mathcal{F}}$, $n\in \mathbb{N}$ be such that
$\lim_{n\to \infty} D(x_{n},x)=0$. From the definition of the
metric $D$, this is obviously equivalent to
$$\lim_{n\to \infty} ||x_{n}^{-} - x_{-}||=
\lim_{n\to \infty} ||x_{n}^{+} - x_{+}||=0,$$ where $||\cdot||$
denotes the uniform norm on $\overline{C}[0,1]$ and
$x_{n}^{-}=(x_{n})_{-}$, $x_{n}^{+}=(x_{n})_{+}$. By the
definition of $x^{*}$, since $L$ is linear  we can write
$$x^{*}(x_{n})=L(x_{n}^{-}) - L(-x_{n}^{+}),$$
where obviously $x_{n}^{-}, -x_{n}^{+}\in IC[0,1], \forall n\in
\mathbb{N}$. Passing to limit with $n\to \infty$ and taking into
account that by hypothesis $L_{|IC[0,1]}$ is continuous with
respect to the uniform convergence on $IC[0,1]$, we immediately
get $x^{*}(x_{n})\to L(x_{-}) - L(-x_{+})=L(x_{-} +
x_{+})=x^{*}(x)$.

Conversely, let $x^{*} : \mathbf{R}_{\mathcal{F}}\to \mathbf{R}$
be linear and continuous on $\mathbf{R}_{\mathcal{F}}$.

Let us consider the set
$$A=\{u\in \overline{C}[0,1] ; \text{ there exist} f, g\in
\overline{C}[0,1], f\le g, f\text{ is increasing on } [0,1],$$
$$g\text{ is decreasing on } [0,1], \text{ such that } u=f+g\}.$$
The set $A$ is a linear subspace of $\overline{C}[0,1]$. Indeed,
for $u, v\in A, u=f+g, v=h+l$ we have $u+v=(f+h)+(g+l)$, where by
hypothesis we easily obtain $f+h\le g+l$, $f+h\in
\overline{C}[0,1]$, $f+h$ is increasing on $[0,1]$, $g+l\in
\overline{C}[0,1]$, $g+l$ is decreasing on $[0,1]$, which implies
that $u+v\in A$. For $\alpha\in \mathbf{R}$ and $u=f+g\in A$, we
have $\alpha u=\alpha f +\alpha g$, where $\alpha f, \alpha g\in
\overline{C}[0,1]$. If $\alpha\ge 0$ then $\alpha f$ is increasing
and $\alpha g$ is decreasing, while if $\alpha < 0$ then $\alpha
g$ is increasing and $\alpha f$ is decreasing but as a
consequence, in both cases it follows $\alpha u\in A$.

First we observe that $IC[0,1]\subset A$. Indeed, for $u\in
IC[0,1]$, we have two possibilities : a) $u(1)\le 0$ or b)
$u(1)>0$. In the case a) we can write $u=u+0\in A$, while in the
case b) we can write $u=f+g$, where $f(t)=u(t)-u(1), g(t)=u(1),
\forall t\in [0,1]$, which proves that $u\in A$. Also, if we
denote by $DC[0,1]$ the set of all $u\in \overline{C}[0,1]$ which
are decreasing on $[0,1]$, then $DC[0,1]\subset A$. Indeed, for
$u\in DC[0,1]$ we have two possibilities : a) $u(1)\ge 0$ when we
write $u=f+g$ with $f=0$, $g=u$ and b) $u(1)<0$ when we write
$u=f+g$ with $f(t)=u(1), \forall t\in [0,1]$ and $g(t)=u(t)-u(1),
\forall t\in [0,1]$.

Now, let us define $L_{0} : A\to \mathbf{R}$ by
$L_{0}(u)=x^{*}(x), \forall u=f+g\in A$, where $x\in
\mathbf{R}_{\mathcal{F}}$ is the unique fuzzy number existing by
Theorem 2.2, such that $x_{-}=f$ and $x_{+}=g$. Notice that if
$u=f+g\in A$, then for all $\varepsilon\ge 0$ we also have the
representation $u=(f-\varepsilon)+(g+\varepsilon)$, where
obviously $f-\varepsilon\le g+\varepsilon$. It follows that for
given $u\in A$, we can choose an infinity of such $x\in
\mathbf{R}_{\mathcal{F}}$, which means that in fact we can define
an infinity of mappings $L_{0}$ as above. For our purposes, we
choose only one, intimately connected to the chosen
representations of the elements $u\in A$, such that for $u\in
IC[0,1]$ and $u\in DC[0,1]$ we choose the above representations.

First we show that $L_{0}$ is linear on $A$. Indeed, for $u=f+g\in
A$, $v=h+l\in A$, we have $L_{0}(u)=x^{*}(x)$,
$L_{0}(v)=x^{*}(y)$, where $x_{-}=f, x_{+}=g$ and $y_{-}=h,
y_{+}=l$. But $(x+y)_{-}=x_{-}+y_{-}=f+h,
(x+y)_{+}=x_{+}+y_{+}=g+l$, which by the linearity of $x^{*}$
implies
$$L_{0}(u+v)=x^{*}(x+y)=x^{*}(x)+x^{*}(y)=L_{0}(u)+L_{0}(v).$$

Now, let $\alpha\in \mathbf{R}$, $u=f+g\in A$ and $x\in
\mathbf{R}_{\mathcal{F}}$ with $x_{-}=f, x_{+}=g$, i.e.
$L_{0}(u)=x^{*}(x)$. If $\alpha\ge 0$ then $(\alpha x)_{-}=\alpha
(x)_{-}=\alpha f$, $(\alpha x)_{+}=\alpha (x)_{+}=\alpha g$, so
$L_{0}(\alpha u)=x^{*}(\alpha x)=\alpha x^{*}(x)=\alpha L_{0}(u)$.
If $\alpha <0$ then $\alpha u=\alpha f + \alpha g=\alpha g +
\alpha f$, where $(\alpha x)_{-}=\alpha g$, $(\alpha x)_{+}=\alpha
f$, which again implies $L_{0}(\alpha u)=\alpha L_{0}(u)$.

As a conclusion, $L_{0}$ is linear on $A$ and by a well-known
result in Functional Analysis (see e.g. [11, pp. 56-57,
Proposition 1.1]), $L_{0}$ can be prolonged to a linear functional
$L : \overline{C}[0,1]\to \mathbf{R}$.

It remains to prove that the restriction of $L_{0}$ to $IC[0,1]$
is continuous on $IC[0,1]$. Thus, let $u_{n}, u\in IC[0,1], n\in
\mathbb{N}$ be such that $u_{n}\to u$, uniformly on $[0,1]$. We
have three possibilities : a) $u(1)<0$ ; b) $u(1)>0$ ; c)
$u(1)=0$.

Case a). We get $u(t)\le u(1)<0, \forall t\in [0,1]$ and by
Theorem 2.2, there is a unique $x\in \mathbf{R}_{\mathcal{F}}$
with $x_{-}(t)=u(t), x_{+}(t)=0, \forall t\in [0,1]$. From
$u_{n}\to u$ uniformly on $[0,1]$, for $\varepsilon
=-\frac{u(1)}{2}>0$, there is $n_{0}\in \mathbb{N}$, such that
$u_{n}(t)-u(t)<-\frac{u(1)}{2}, \forall t\in [0,1], n\ge n_{0}$,
which implies that for all $t\in [0,1]$ and $n\ge n_{0}$ we have
$$u_{n}(t)<u(t)-\frac{u(1)}{2}\le
u(1)-\frac{u(1)}{2}=\frac{u(1)}{2}\le 0.$$ Therefore, for any
$n\ge n_{0}$ there is $x_{n}\in \mathbf{R}_{\mathcal{F}}$ such
that $x_{n}^{-}(t)=u_{n}(t), x_{n}^{+}(t)=0, \forall t\in [0,1]$.
As a conclusion, $u_{n}\to u$ uniformly on $[0,1]$, implies
$x_{n}^{-}\to x_{-}$ and $x_{n}^{+}\to x_{+}$, uniformly on
$[0,1]$, which by e.g. [13, p. 524] implies $D(x_{n},x)\to 0$
(when $n\to \infty$) and together with the continuity of $x^{*}$
we get
$$L_{0}(u_{n})=x^{*}(x_{n})\to x^{*}(x)=L_{0}(u).$$

Case b). There is a unique $x\in \mathbf{R}_{\mathcal{F}}$ with
$x_{-}(t)=u(t)-u(1), x_{+}(t)=u(1), \forall t\in [0,1]$. Since
$\lim_{n\to \infty}u_{n}(1)=u(1)>0$, there exists $n_{0}\in
\mathbf{N}$ such that $u_{n}(1)>0, \forall n\ge n_{0}$. Therefore,
for any $n\ge n_{0}$, there is a unique $x_{n}\in
\mathbf{R}_{\mathcal{F}}$ with $x_{n}^{-}(t)=u_{n}(t)-u_{n}(1),
x_{n}^{+}(t)=u_{n}(1), \forall t\in [0,1]$. Then, obviously again
we obtain $x_{n}^{-}(t)\to u(t)-u(1)=x_{-}$ and $x_{n}^{+}\to
u(1)=x_{+}$, uniformly on $[0,1]$, which implies $D(x_{n},x)\to 0$
(when $n\to \infty$) and by the continuity of $x^{*}$ it follows
$$L_{0}(u_{n})=x^{*}(x_{n})\to x^{*}(x)=L_{0}(u).$$

Case c). From $u(t)\le u(1)=0, \forall t\in [0,1]$, there exists a
unique $x\in \mathbf{R}_{\mathcal{F}}$ with $x_{-}(t)=u(t),
x_{+}(t)=0, \forall t\in [0,1]$. Concerning each term of the
sequence $(u_{n})_{n}$, we have two possibilities :

(i) $u_{n}(1)\le 0$ or (ii) $u_{n}(1)>0$.

Subcase (i). We have $u_{n}(t)\le 0, \forall t\in [0,1]$ and there
is a unique $x_{n}\in \mathbf{R}_{\mathcal{F}}$ with
$x_{n}^{-}(t)=u_{n}(t), x_{n}^{+}(t)=0, \forall t\in [0,1]$.

Subcase (ii). There is a unique $x_{n}\in
\mathbf{R}_{\mathcal{F}}$ with $x_{n}^{-}(t)=u_{n}(t)-u_{n}(1),
x_{n}^{+}(t)=u_{n}(1), \forall t\in [0,1]$.

From both subcases we obtain that $u_{n}(t)\to u(t)$ uniformly on
$[0,1]$ implies $x_{n}^{-}(t)\to u(t)=x_{-}(t), x_{n}^{+}(t)\to
u(1)=0=x_{+}(t)$, uniformly on $[0,1]$, i.e. $\lim_{n\to
\infty}D(x_{n},x)=0$ and reasoning as for the above cases we
obtain the continuity of the restriction of $L_{0}$ to $IC[0,1]$
in this last subcase too. The theorem is proved.

{\bf REMARK}. A natural family of functionals $L$ in the statement
of Theorem 3.1 can be defined as follows. For any fixed continuous
function $h : [0,1]\to \mathbf{R}$, first define $L_{0} : A\to
\mathbf{R}$, as the Riemann-Stieltjes integral
$L_{0}(u)=\int_{0}^{1}h(t)d[u(t)]$. Then $L$ will be a linear
extension of $L_{0}$ to $\overline{C}[0,1]$. Obviously that for
any $u=f+g\in A$, $L_{0}(u)$ has sense and we have
$L_{0}(u)=\int_{0}^{1}h(t)d[f(t)]+\int_{0}^{1}h(t)d[g(t)]$.

It is easy to prove that $L_{0}$ is linear on $A$. Also, the
restriction of $L_{0}$ to $IC[0,1]$ is continuous on $IC[0,1]$.
Indeed, let $u_{n}, u\in IC[0,1], n\in \mathbb{N}$, $u_{n}\to u$
uniformly on $[0,1]$ (when $n\to \infty$). We have : $\bigvee
_{0}^{1}u_{n}=u_{n}(1)-u_{n}(0)\to u(1)-u(0)$, when $n\to \infty$.
This means that there exists $M>0$ such that $\bigvee
_{0}^{1}u_{n}\le M, \forall n\in \mathbb{N}$, which by the
classical Helly-Bray theorem (see e.g. [12, p. 38]) implies
$$\lim_{n\to \infty}L_{0}(u_{n})=\lim_{n\to
\infty}\int_{0}^{1}h(t)d[u_{n}(t)]=\int_{0}^{1}h(t)d[u(t)]=L_{0}(u),$$
proving the desired continuity. Therefore, a class of linear
continuous functionals $x^{*} : \mathbf{R}_{\mathcal{F}}\to
\mathbf{R}$ are of the form
$$x^{*}(x)=\int_{0}^{1}h(t)d[x_{-}(t)]+\int_{0}^{1}h(t)d[x_{+}(t)],
\forall x\in \mathbf{R}_{\mathcal{F}},$$ where $h : [0,1]\to
\mathbf{R}$ is continuous on $[0,1]$.

Notice that from a well-known formula for the Riemann-Stieltjes
integral (see e.g. [12, p. 30]), we also can write
$$x^{*}(x)=h(1)[x_{-}(1)+x_{+}(1)]-h(0)[x_{-}(0)+x_{+}(0)]
-\int_{0}^{1}[x_{-}(t)+x_{+}(t)]d[h(t)].$$ It is natural the
following.

{\bf OPEN QUESTION 1}. Are all the linear continuous functionals
$x^{*} : \mathbf{R}_{\mathcal{F}}\to \mathbf{R}$ of the form in
the previous remark ?

{\bf REMARK.} A method to answer the above open question would be
to use the ideas in the proof of classical Riesz's result
concerning the form of linear continuous functionals on $C[0,1]$.
Unfortunately, it seems that these ideas do not work in our case,
since if we define $h(t)=L_{0}(z_{t}), t\in [0,1]$, (where
$z_{t}(x)=1, x\in [0,t), z_{t}(x)=0, x\in [t,1]$ and the
restriction  of $L_{0}$ to $IC[0,1]$ is supposed to be continuous
on $IC[0,1]$ with respect to the uniform convergence), then $h$ is
not, in general, continuous on $[0,1]$.

In order to prove our second main result, we need the following
notations:
$$SIC[0,1]=\{(s_{n})_{n} ; s_{n}\in IC[0,1], \forall n\in
\mathbb{N}, (s_{n})_{n} \text { is uniformly convergent }\},$$
$$S\overline{C}[0,1]=\{(s_{n})_{n} ; s_{n}\in \overline{C}[0,1],
\forall n\in \mathbb{N}, (s_{n})_{n} \text { is uniformly
convergent }\}.$$

{\bf THEOREM 3.2 } {\it $ x^{*} : c_{\mathbf{R}_{\mathcal{F}}}\to
\mathbf{R}$ is a linear continuous functional on
$c_{\mathbf{R}_{\mathcal{F}}}$, if and only if there exists a
linear functional $L : S\overline{C}[0,1]\to \mathbf{R}$ such that
$$x^{*}(x)=L(x_{-}+x_{+}), \forall x=(x_{n})_{n}\in
c_{\mathbf{R}_{\mathcal{F}}},$$ where the restriction of $L$ to
$SIC[0,1]$ is continuous on $SIC[0,1]$ with respect to the
convergence induced by the metric on $S\overline{C}[0,1]$ defined
by
$$\Phi[(s_{n})_{n},(t_{n})_{n}]=\sup\{||s_{n} - t_{n}|| ; n\in
\mathbb{N}\}.$$ Here $||\cdot||$ represents the uniform norm, for
$x=(x_{n})_{n}\in c_{\mathbf{R}_{\mathcal{F}}}$, we have denoted
$x_{-}=(x_{n}^{-})_{n}$, $x_{+}=(x_{n}^{+})_{n}$ and obviously
that the convergence of the sequence $x=(x_{n})_{n}$ in the metric
$\mu$ in $c_{\mathbf{R}_{\mathcal{F}}}$, implies that $x_{-}$ and
$x_{+}$ are uniformly convergent sequences of functions.}

{\it Proof}. First, let us suppose that $x^{*} :
c_{\mathbf{R}_{\mathcal{F}}}\to \mathbf{R}$ is of the form in
statement. The linearity of $x^{*}$ follows exactly as in the
proof of Theorem 3.1. To prove the continuity of $x^{*}$, let
$x=(x_{n})_{n}\in c_{\mathbf{R}_{\mathcal{F}}}$,
$x^{m}=(x_{n}^{(m)})_{n}\in c_{\mathbf{R}_{\mathcal{F}}},
m=1,2,...,$ be such that $\lim_{m\to \infty}\mu(x^{m},x)=0$. From
the definition of $\mu$ (see Section 2) and $\Phi$, it is
immediate that
$$\Phi[(x^{m})_{-},x_{-}]\to 0, \Phi[(x^{m})_{+},x_{+}]\to 0,$$
when $m\to \infty$. Since we can write
$x^{*}(x)=L(x_{-})-L(-x_{+})$ and by $(x^{m})_{-}, x_{-},
-(x^{m})_{+}, -x_{+}\in SIC[0,1]$, obviously that the continuity
of $L$ implies the continuity of $x^{*}$. (Above we have denoted
$(x^{m})_{-}=((x_{n}^{m})_{-})_{n},
(x^{m})_{+}=((x_{n}^{m})_{+})_{n}, \forall m=1,2,...$).

Conversely, let $x^{*} : c_{\mathbf{R}_{\mathcal{F}}}\to
\mathbf{R}$ be linear continuous functional on
$c_{\mathbf{R}_{\mathcal{F}}}$. We use a similar idea to that in
the proof of Theorem 3.1. Firstly, it is easy to show by standard
procedure that with respect to usual addition and scalar
multiplication of the sequences of functions and the norm
$||(s_{n})_{n}||=\sup\{||s_{n}(t)|| ; n\in \mathbb{N}\}$,
$S\overline{C}[0,1]$ becomes a real Banach space. Also, we have
$\Phi[(s_{n})_{n},(t_{n})_{n}]=||(s_{n})_{n} - (t_{n})_{n}||$.

Now, let us define the set
$$SA=\{u=(u_{n})_{n} ; u_{n}=f_{n}+g_{n}\in A, \text { such that }
f=(f_{n})_{n}, g=(g_{n})_{n}$$
$$\text { are uniformly convergent }\},$$
where the set $A$ is defined in the proof of Theorem 3.1.
As in the previous proof, we easily obtain that $SA$ is a linear
subspace of $S\overline{C}[0,1]$. We define $L_{0} : SA \to
\mathbf{R}$ by $L_{0}(u)=x^{*}(x), \forall u=f+g\in SA$, where
$x=(x_{n})_{n}\in c_{\mathbf{R}_{\mathcal{F}}}$ is the unique
sequence of fuzzy numbers satisfying the relations
$(x_{n})_{-}=f_{n}, (x_{n})_{+}=g_{n}, \forall n\in \mathbb{N}$.

The linearity of $L_{0}$ follows exactly as in the proof of
Theorem 3.1 and therefore there exists a linear extension of
$L_{0}$, denoted by $L : S\overline{C}[0,1]\to \mathbf{R}$.

On the other hand, from the inclusion $IC[0,1]\subset A$, we
immediately get the inclusion $SIC[0,1]\subset SA$. It remains to
prove the continuity on $SIC[0,1]$ of the restriction of $L_{0}$
to $SIC[0,1]$. Thus, let $u^{m}, u\in SIC[0,1]$ be such that
$\lim_{m\to \infty}\Phi(u^{m},u)=0$. Denoting
$u^{m}=(u^{m}_{n})_{n}, u=(u_{n})_{n}$, this means
$u_{n}^{m}(t)\to u_{n}(t)$ (with $m\to \infty$), uniformly with
respect to $t\in [0,1]$ and $n\in \mathbb{N}$. According to the
proof of Theorem 3.1, there exists (in a unique way) $x_{n}\in
c_{\mathbf{R}_{\mathcal{F}}}$ (depending on $u_{n}$) and
$x_{n}^{m}\in c_{\mathbf{R}_{\mathcal{F}}}$ (depending on
$u_{n}^{m}$, such that $D(x_{n}^{m},x_{n})\to 0$, (with $m\to
\infty$), uniformly respect to $n\in \mathbb{N}$, where
$x=(x_{n})_{n}, x^{m}=(x_{n}^{m})_{n}\in
c_{\mathbf{R}_{\mathcal{F}}}$, for all $m=1,2,...,$ and
$$L_{0}(u^{m})=x^{*}(x^{m}), x^{*}(x)=L_{0}(u).$$
From the continuity of $x^{*}$ it follows the continuity of
$L_{0}$. Note that as in the proof of Theorem 3.1, we define here
$L_{0}$ under the hypothesis that for $u\in IC[0,1]$ and $u\in
DC[0,1]$ we choose the representations in the proof of Theorem
3.1, which implies the corresponding representation for the
elements in $SIC[0,1]$.

{\bf REMARK.} A natural class of linear continuous functionals
$x^{*} : c_{\mathbf{R}_{\mathcal{F}}}\to \mathbf{R}$ is given by
the form
$$x^{*}(x)=\int_{0}^{1}h_{1}(t)d[z_{-}(t)+z_{+}(t)]+
\sum_{j=1}^{\infty}\alpha_{j}\int_{0}^{1}h_{2}(t)d[(x_{j})_{-}(t)
+(x_{j})_{+}(t)],$$ for all $x=(x_{j})_{j}\in
c_{\mathbf{R}_{\mathcal{F}}}$ with $\lim_{j\to
\infty}D(x_{j},z)=0$, where $h_{1}, h_{2} : [0,1]\to \mathbf{R}$
are arbitrary continuous functions on $[0,1]$ and
$(\alpha_{j})_{j}$ is an arbitrary sequence of real numbers
satisfying $\sum_{j=1}^{\infty}|\alpha_{j}|<+\infty$.

Taking into account the well-known "inversion" formula for the
Riemann-Stieltjes integral, $x^{*}(x)$ in the above formula can be
written by
$$x^{*}(x)=h_{1}[z_{-}(1)+z_{+}(1)]-h_{1}(0)[z_{-}(0)+z_{+}(0)]-
\int_{0}^{1}[z_{-}(t)+z_{+}(t)]d[h_{1}(t)]+$$
$$\sum_{j=1}^{\infty}\alpha_{j}\{h_{2}(1)[(x_{j})_{-}(1)+(x_{j})_{+}(1))-
h_{2}(0)((x_{j})_{-}(0)+(x_{j})_{+}(0)]-$$
$$\int_{0}^{1}[(x_{j})_{-}(t)+(x_{j})_{+}(t)]d[h_{2}(t)]\}.$$

{\bf OPEN QUESTION 2 } It is an open question if all the linear
continuous functionals on $c_{\mathbf{R}_{\mathcal{F}}}$ are of
the above form.

In a similar manner can be proved the following four theorems.

{\bf THEOREM 3.3 } {\it $ x^{*} : m_{\mathbf{R}_{\mathcal{F}}}\to
\mathbf{R}$ is a linear continuous functional, if and only if
there exists a linear functional $L : M\overline{C}[0,1]\to
\mathbf{R}$ such that
$$x^{*}(x)=L(x_{-}+x_{+}), \forall x=(x_{n})_{n}\in
m_{\mathbf{R}_{\mathcal{F}}},$$ where the restriction to
$MIC[0,1]$ of $L$ is continuous with respect to the metric $\Phi$
in Theorem 3.2.

Here
$$M\overline{C}[0,1]=\{(s_{n})_{n} ; s_{n}\in
\overline{C}[0,1], \forall n\in \mathbb{N}, (s_{n})_{n} \text {is
uniformly bounded }\},$$
$$MIC[0,1]=\{s_{n}\in IC[0,1], \forall n\in \mathbb{N}, (s_{n})_{n} \text
{is uniformly bounded }\},$$ and for $x=(x_{n})_{n}\in
m_{\mathbf{R}_{\mathcal{F}}}$, we have denoted
$x_{-}=((x_{n})_{-})_{n}, x_{+}=((x_{n})_{+})_{n}$, which also are
uniformly bounded.}

{\bf THEOREM 3.4 } {\it Let $1\le p<+\infty$. Then $x^{*} :
l^{p}_{\mathbf{R}_{\mathcal{F}}}\to \mathbf{R}$ is linear and
continuous functional, if and only if there exists a linear
functional $L : S^{p}\overline{C}[0,1]\to \mathbf{R}$ such that
$$x^{*}(x)=L(x_{-}+x_{+}), \forall x=(x_{n})_{n}\in
l^{p}_{\mathbf{R}_{\mathcal{F}}},$$ where the restriction to
$S^{p}IC[0,1]$ of $L$ is continuous with respect to the metric on
$S^{p}\overline{C}[0,1]$
$$\Psi[(s_{n})_{n},(t_{n})_{n}]=\left\{\sum_{n=1}^{\infty}||s_{n} -
t_{n}||^{p}\right\}^{1/p}.$$ Here $||\cdot||$ denotes the uniform
norm, for $x=(x_{n})_{n}\in l^{p}_{\mathbf{R}_{\mathcal{F}}}$ the
sequences of functions $x_{-}=((x_{n})_{-})_{n},
x_{+}=((x_{n})_{+})_{n}$ satisfy
$$\sum_{n=1}^{\infty}||(x_{n})_{-}||^{p}<+\infty,
\sum_{n=1}^{\infty}||(x_{n})_{+}||^{p}<+\infty,$$ and we have the
notations
$$S^{p}\overline{C}[0,1]=\{(s_{n})_{n} ; s_{n}\in
\overline{C}[0,1], \sum_{n=1}^{\infty}||s_{n}||^{p}<+\infty\},$$
$$S^{p}IC[0,1]=\{(s_{n})_{n} ; s_{n}\in
IC[0,1], \sum_{n=1}^{\infty}||s_{n}||^{p}<+\infty\}.$$}

{\bf REMARK.} A class of linear continuous functionals
$x^{*}:l^{p}_{\mathbf{R}_{\mathcal{F}}}\to \mathbf{R}$ is given by
the formula
$$x^{*}(x)=\sum_{j=1}^{\infty}\alpha_{j}\int_{0}^{1}h(t)
d[(x_{j})_{-}(t)+(x_{j})_{+}(t)], \forall x=(x_{j})_{j}\in
l^{p}_{\mathbf{R}_{\mathcal{F}}},$$ where $h : [0,1]\to
\mathbf{R}$ is continuous on $[0,1]$ and $(\alpha_{j})_{j}$ is a
sequence of real numbers satisfying : (i) $|\alpha_{j}|\le M,
\forall j\in \mathbb{N}$ if $p=1$ and (ii) $\sum_{j=1}^{\infty}
|\alpha_{j}|^{q}<+\infty, 1/p + 1/q = 1$, if $1<p<+\infty$.

{\bf OPEN QUESTION 3 } It is an open question if the linear
continuous functionals on $l^{p}_{\mathbf{R}_{\mathcal{F}}}$ are
all of the above form.

{\bf THEOREM 3.5 } {\it $ x^{*} :
C([a,b];\mathbf{R}_{\mathcal{F}})\to \mathbf{R}$ is linear
continuous functional if and only if there exists a linear
functional $L : C([a,b];\overline{C}[0,1])\to \mathbf{R}$, such
that
$$x^{*}(x)=L(x_{-}+x_{+}), \forall x\in
C([a,b];\mathbf{R}_{\mathcal{F}}),$$ where the restriction to
$CIC[0,1]$ of $L$ is continuous on $CIC[0,1]$ with respect to the
metric on $C([a,b];\overline{C}[0,1])$ given by
$$\Delta(F,G)=\sup\{||F(t) - G(t)|| ; t\in [a,b]\},$$
for all $F, G : [a,b]\to \overline{C}[0,1]$ continuous on $[a.b]$,
where $||\cdot||$ is the uniform norm on $\overline{C}[0,1]$ and
$\overline{C}[0,1]$ is considered endowed with the uniform metric
$\Gamma(f,g)=||f - g||$.

Here
$$C([a,b];\overline{C}[0,1])=\{F : [a,b]\to \overline{C}[0,1] ; F
\text { is continuous on } [a,b]\},$$
$$CIC[0,1]=\{F\in C([a,b];\overline{C}[0,1]) ; F(t)\in IC[0,1],
\forall t\in [a,b]\},$$ and for $x\in C([a,b] ;
\mathbf{R}_{\mathcal{F}})$, we define $x_{-}, x_{+} : [a,b]\to
\overline{C}[0,1]$ by
$$[x_{-}(t)](r)=[x(t)]_{-}(r), [x_{+}(t)](r)=[x(t)]_{+}(r),
\forall t\in [a,b], r\in [0,1],$$ which obviously satisfy $x_{-},
x_{+} \in C([a,b] ; \overline{C}[0,1])$.}

{\bf REMARK.} A class of linear continuous functionals $x^{*} :
C([a,b];\mathbf{R}_{\mathcal{F}})\to \mathbf{R}$ is given by the
formula
$$x^{*}(x)=\int_{a}^{b}\left\{\int_{0}^{1}h_{1}(s)d[x(t)_{-}(s) +
x(t)_{+}(s)]\right\}d[h_{2}(t)],$$ where $h_{1} ; [0,1]\to
\mathbf{R}$ is continuous on $[0,1]$ and $h_{2} ; [a,b]\to
\mathbf{R}$ is of bounded variation, arbitrary.

{\bf OPEN QUESTION 4 } Remains an open question if all the linear
continuous functionals on $C([a,b];\mathbf{R}_{\mathcal{F}})$ are
of this form.

{\bf THEOREM 3.6 } {\it Let $1\le p< +\infty$. Then $x^{*} :
L^{p}([a,b] ; \mathbf{R}_{\mathcal{F}})\to \mathbf{R}$ is linear
continuous functional, if and only if there exists a linear
functional $L : L^{p}([a,b] ; \overline{C}[0,1])\to \mathbf{R}$,
such that
$$x^{*}(x)=L(x_{-}+x_{+}),$$
where the restriction to $L^{p}IC[0,1]$ of $L$ is continuous on
$L^{p}IC[0,1]$ with respect to the metric on
$L^{p}([a,b];\overline{C}[0,1])$ given by
$$\Delta_{p}(F,G)=\left\{\int_{a}^{b}||F(t) - G(t)||^{p}\right\}^{1/p},$$
for all $F, G \in L^{p}([a,b];\overline{C}[0,1])$, with
$||\cdot||$ the uniform norm on $\overline{C}[0,1]$ and
$\int_{a}^{b}$ the Lebesgue-kind integral.

Here
$$L^{p}([a,b];\overline{C}[0,1])=\{F : [a,b]\to \overline{C}[0,1]
; \int_{a}^{b}||F(t)||^{p}dt < +\infty\},$$
$$L^{p}IC[0,1]=\{F\in L^{p}([a,b];\overline{C}[0,1]) ; F(t)\in
IC[0,1], \forall t\in [0,1]\},$$ and $x_{-}, x_{+}$ are defined as
in the statement of Theorem 3.5.}

{\bf REMARKS.} 1) A class of linear continuous functionals $x^{*}
: L^{p}([a,b] ; \mathbf{R}_{\mathcal{F}})\to \mathbf{R}$ is given
by the formula
$$x^{*}(x)=\int_{a}^{b}\left\{\int_{0}^{1}h_{1}(s)d[x(t)_{-}(s)+
x(t)_{+}(s)]\right\}h_{2}(t)dt,$$ where $h_{1} : [0,1]\to
\mathbf{R}$ is continuous on $[0,1]$, $h_{2}\in
L^{q}([a,b];\mathbf{R})$ with $1/p + 1/q=1$ if $1<p<+\infty$ and
$h_{2}$ is a.e. bounded on $[a,b]$, in the case when $p=1$.

{\bf OPEN QUESTION 5 } Remains an open question if all the linear
continuous functionals on $L^{p}([a,b] ;
\mathbf{R}_{\mathcal{F}})$ are of the form in Remark 1.

2) A crucial step in the proofs of Theorems 3.1 and 3.2 is the
construction of the set $A$ and the linear functional $L_{0} :
A\to \mathbf{R}$ and of $SA$ and $L_{0} : SA\to \mathbf{R}$,
respectively.

In the case of Theorem 3.5, the corresponding constructions are
the set
$$FA=\{U\in C([a,b];\mathbf{R}_{\mathcal{F}}) ; \text { there
exist } F,G \text { with } F, -G\in CIC[0,1] \text { such that }$$
$$F(t)(r)\le G(t)(r), \forall t\in [a,b], r\in [0,1] \text { and }
U(t)=F(t) + G(t), \forall t\in [a,b]\},$$ and $L_{0} : FA\to
\mathbf{R}$ is defined by $L_{0}(U)=x^{*}(V)$, where $x^{*} :
C([a,b];\mathbf{R}_{\mathcal{F}})\to \mathbf{R}$ is given and $V$
is the unique function $V : [a,b]\to \mathbf{R}_{\mathcal{F}}$
obtained (by Theorem 2.2) from the relations
$$[V(t)]_{-}=F(t), [V(t)]_{+}=G(t), \forall t\in [a,b].$$

In the case of Theorem 3.6, the construction of $FA$ and $L_{0}$
is similar, with the difference that $U\in
L^{p}([a,b];\mathbf{R}_{\mathcal{F}})$ and $x^{*} :
L^{p}([a,b];\mathbf{R}_{\mathcal{F}})\to \mathbf{R}$.

Of course that in both cases (of Theorems 3.5 and 3.6), the
mapping $L_{0}$ is defined under the hypothesis that for $u\in
IC[0,1]$ (and $u\in DC[0,1]$) we adopt the representations in the
proof of Theorem 3.1.

\end{document}